\documentclass{amsart}
\usepackage{graphicx}
\usepackage{calc}
\usepackage{hyperref}

\usepackage[cmtip,arrow,all,2cell]{xy}
\usepackage{pb-diagram,pb-xy}
\usepackage{amsmath,amsthm,amssymb,tabularx}

\theoremstyle{plain}

\newtheorem{tm}{Theorem}[section]

\newtheorem{prop}{Proposition}[section]

\newtheorem{lm}{Lemma}[section]

\newtheorem{cor}{Corollary}[section]

\theoremstyle{definition}

\newtheorem{ex}{Example}[section]

\newtheorem{rem}[ex]{Remark}

\newcommand{\beq}{\begin{equation}}
\newcommand{\eeq}{\end{equation}}
\newcommand{\bit}{\begin{itemize}}
\newcommand{\eit}{\end{itemize}}
\newcommand{\btm}{\begin{tm}}
\newcommand{\etm}{\end{tm}}
\newcommand{\blm}{\begin{lm}}
\newcommand{\elm}{\end{lm}}
\newcommand{\bcor}{\begin{cor}}
\newcommand{\ecor}{\end{cor}}
\newcommand{\bex}{\begin{ex}}
\newcommand{\eex}{\end{ex}}
\newcommand{\bcx}{\begin{cex}}
\newcommand{\ecx}{\end{cex}}
\newcommand{\bers}{\begin{ers}}
\newcommand{\eers}{\end{ers}}
\newcommand{\ber}{\begin{er}}
\newcommand{\eer}{\end{er}}
\newcommand{\bdf}{\begin{df}}
\newcommand{\edf}{\end{df}}
\newcommand{\brem}{\begin{rem}}
\newcommand{\erem}{\end{rem}}
\newcommand{\bpr}{\begin{proof}}
\newcommand{\epr}{\end{proof}}
\newcommand{\bprop}{\begin{prop}}
\newcommand{\eprop}{\end{prop}}

\newcommand{\norm}[1]{\left\Vert#1\right\Vert}

\def \N {\mathbb{N}}

\def \R {\mathbb{R}}

\def \C {\mathbb{C}}

\def\Ste{\operatorname{\sf{Ste}}}

\def \e {\varepsilon}
\def \ph {\varphi}

\def\Ob{\operatorname{\sf Ob}}
\def\id{\operatorname{\sf id}}
\def\Mor{\operatorname{\sf Mor}}

\def\DEpi{\operatorname{\sf DEpi}}

\def\env{\operatorname{\sf env}}

\def\Env{\operatorname{\sf Env}}
\def\Rf{\operatorname{\sf Ref}}

\def\Im{\operatorname{\sf Im}}
\def\ker{\operatorname{\sf ker}}
\def\Ker{\operatorname{\sf Ker}}
\def\coker{\operatorname{\sf coker}}
\def\Coker{\operatorname{\sf Coker}}
\def\Dom{\operatorname{\sf Dom}}
\def\Ran{\operatorname{\sf Ran}}
\def\Aug{\operatorname{\sf Aug}}

\def\AugC*{\operatorname{\sf AugC^*}}

\def\InvSteAlg{\operatorname{\sf InvSteAlg}}
\def\AugInvSteAlg{\operatorname{\sf AugInvSteAlg}}
\def\Arr{\operatorname{\sf Arr}}

\begin{document}

\title[Kernels and cokernels in the category of stereotype algebras]{Kernels and cokernels in the category of augmented involutive stereotype algebras}

\author{S.S.Akbarov}
 \address{School of Applied Mathematics, National Research University Higher School of Economics, 34, Tallinskaya St. Moscow, 123458 Russia}
 \email{sergei.akbarov@gmail.com}

\thanks{Supported by the RFBR grant No. 18-01-00398.}

\maketitle

\begin{abstract}
We prove several properties of kernels and cokernels in the category of augmented involutive stereotype algebras: 1) this category has kernels and cokernels, 2) the cokernel is preserved under the passage to the group stereotype algebras, and 3) the notion of cokernel allows to prove that the continuous envelope $\Env{\mathcal C}^\star(Z\cdot K)$ of the group algebra of a compact buildup of an abelian locally compact group is an involutive Hopf algebra in the category of stereotype spaces $(\Ste,\odot)$. The last result plays an important role in the generalization of the Pontryagin duality for arbitrary Moore groups.
\end{abstract}

\section{Introduction}

The Pontryagin duality theorem for locally compact commutative groups has been generalized many times to various classes of topological groups, not necessarily commutative \cite{Akbarov-Shavgulidze,Chasco-Dikranjan-Martin-Peinador,Enock-Schwartz,Hernandez-Arrieta,Kustermans-Vaes,Timmermann}.

The first generalizations to commutative groups appeared in the works by S.~Kaplan of 1948 and 1950 \cite{Kaplan-1948,Kaplan-1950}, and by now it is known in particular, that wide classes of topological vector spaces can be considered as examples of groups on which this duality acts \cite{Akbarov-Shavgulidze}. One of those classes is the class of stereotype spaces described by the author in \cite{Akbarov}.

The first generalization to the groups without the condition of commutativity was suggested in the works by T.~Tannaka \cite{Tannaka} of 1938 and by M.~G.~Krein \cite{Krein} of 1949. They built a duality theory for arbitrary (not necessarily commutative) compact groups. In 1973 L.~I.~Vainerman, G.~I.~Kac, M.~Enock and J.-M.~Schwartz built a general theory for all locally compact groups (see details in the book \cite{Enock-Schwartz}). From 1980-ies the research in this area was resumed after   discovery of quantum groups, to which the constructed theories began to be actively transferred (see \cite{Kustermans-Vaes,Timmermann}).

To date, the Vainerman---Kac--- Enock---Schwartz theory and its extensions cover a very wide class of (non necessarily commutative) groups and group-like objects, however, as noted in  \cite{Akbarov-stein-groups}, all these generalizations have the disadvantage that the enclosing category (containing the category of locally compact groups) consists of objects that resemble Hopf algebras (and even echo with them in names, such as {\it ``Hopf---von Neumann algebra''} or {\it ``$C^*$-Hopf algebra''}), but formally {\it are not Hopf algebras}.

The desire to find a theory free from this flaw led the author to the investigations in this area as part of a project on the study of applications of stereotype spaces, announced in \cite{Akbarov-3}. As an alternative, in 2008 the author constructed a duality theory for a class of compactly generated complex Lie groups with an affine algebraic connected component of unit \cite{Akbarov-stein-groups}. The main result of this article is described by the diagram
 \beq\label{int:chetyrehugolnik-O-O*}
 \xymatrix @R=2.pc @C=4.pc @M=14pt
 {
 {\mathcal O}^\star(G)
  \ar@{|->}[r]^{\heartsuit} &
 {\mathcal O}_{\exp}^\star(G)
 \ar@{|->}[d]^{\star}
 \\
 {\mathcal O}(G) \ar@{|->}[u]^{\star}
 &
 {\mathcal O}_{\exp}(G) \ar@{|->}[l]_{\heartsuit}
 }
 \eeq
where $G$ is an arbitrary compactly generated complex Lie group with an affine algebraic component of unit, ${\mathcal O}(G)$ the algebra of holomorphic functions on it, ${\mathcal O}^\star(G)$ its dual algebra of analytic functionals, ${\mathcal O}_{\exp}(G)$ the algebra of holomorphic functions of exponential type on $G$, ${\mathcal O}_{\exp}^\star(G)$ its dual algebra of exponential analytic functionals, $\star$ the passage to the dual stereotype space, and $\heartsuit$ the passage to the Arens-Michael envelope.

The key detail in diagram \eqref{int:chetyrehugolnik-O-O*} is that the objects in its upper left and lower right corners --- ${\mathcal O}^\star(G)$ and ${\mathcal O}_{\exp}(G)$ --- are Hopf algebras in a certain monoidal category, namely, the category of stereotype spaces $({\tt Ste},\circledast)$. This allows us to interpret \eqref{int:chetyrehugolnik-O-O*} as a construction describing a generalization of Pontryagin duality to the class of groups in question (see a detailed discussion in \cite{Akbarov-stein-groups}):
\bit{

\item[1)] the fact that the movement in diagram \eqref{int:chetyrehugolnik-O-O*} at the fourth step returns us back to the original place can be naturally interpreted as a duality identity
\beq\label{widehat-widehat-H-cong-H}
\widehat{\widehat{H}}\cong H,
\eeq
for the functor $H\mapsto\widehat{H}$ of composition of operations $\heartsuit$ and $\star$
    $$
    \widehat{H}=(H^\heartsuit)^\star,
    $$
on the Hopf algebras $H={\mathcal O}^\star(G)$ and $H={\mathcal O}_{\exp}(G)$,

\item[2)] the operation $G\mapsto {\mathcal O}^\star(G)$ that assigns to a group its group algebra, can be understood as an embedding of the considered category of groups into the category of Hopf algebras over $({\tt Ste},\circledast)$ with the property \eqref{widehat-widehat-H-cong-H}, and

\item[3)] the fact that the resulting embedding of a category of groups into the category of ``reflexive in the sense of \eqref{widehat-widehat-H-cong-H} Hopf algebras'' generalizes Pontryagin duality is expressed by the identity for Abelian groups $G$ \cite[Theorem 7.3]{Akbarov-stein-groups}
    $$
    {\mathcal O}^\star(\widehat{G})\cong \widehat{{\mathcal O}^\star(G)}
    $$
    (here $\widehat{G}$ means the Pontryagin dual group).
}\eit
As a result, we obtain a generalization of Pontryagin's duality, free from the flaws of the Vainerman–--Kac---Enock---Schwartz theory. Although it acts on a much narrower class of groups, it is still interesting since it is much simpler and {\it it covers all complex affine algebraic groups}.

In 2013, Yu.~N.~Kuznetsova \cite{Kuznetsova} made an attempt to construct a similar theory for the class of arbitrary Moore groups, where the Arens-Michael envelope $\heartsuit$ was replaced by another operation, called {\it $C^*$-envelope}. Later, however, it turned out that Kuznetsova's work contains an error, and its main results can be considered proven only for groups of the form $\R^n\times K \times D$, where $n\in\N$, $K$ is a compact group, and $D$ a discrete Moore group (see \cite[Errata]{Akbarov-C^infty-2}).

This paper is devoted to some preliminary steps for elimination the error in  \cite{Kuznetsova} within the framework of the general scheme of duality described by the author in \cite{Akbarov-C^infty-1,Akbarov-C^infty-2}. We prove here an important statement which is necessary for generalization of Kuznetsova's result to the class of arbitrary Moore groups: a continuous envelope (an analog of the $C^*$-envelope of \cite{Kuznetsova}) $\Env_{\mathcal C}{\mathcal C}^\star(Z\cdot K)$ of a group algebra of the form $Z\cdot K$, where $Z$ is an Abelian locally compact group and $K$ is a compact group, is an involutive Hopf algebra in the category of stereotype spaces $(\Ste,\odot)$. We will need this result in \cite{Akbarov-Moore}, where it will be a main tool in proving the fact that algebras arising in a diagram similar to \eqref{int:chetyrehugolnik-O-O*} are Hopf algebras in the category ${\tt Ste}$.

Our reasonings are based on the properties of kernels and cokernels in the category $\AugInvSteAlg$ of augmented involutive stereotype algebras.

\section{Kernel and cokernel in a category with a zero object}

Let us recall \cite{General-algebra,MacLane}, that a {\it zero object} or {\it zero} in a category ${\tt K}$ is an object $0$ such that for each object $X$ in ${\tt K}$ there exists a unique morphism $X\to 0$, and a unique morphism $0\to X$. If a category ${\tt K}$ has a zero object, then it is unique up to an isomorphism. A morphism  $\ph:X\to Y$ in a category ${\tt K}$ with a zero object $0$ is called a {\it zero morphism}, if it can be factored through the zero object, i.e. if $\ph$ is a composition of two (unique) morphism $X\overset{0_{X,0}}{\longrightarrow}0$ and $0\overset{0_{0,Y}}{\longrightarrow}Y$:
$$
\begin{diagram}
\node{X}\arrow{se,b}{0_{X,0}}\arrow[2]{e,t}{\ph}\node[2]{Y} \\
\node[2]{0}\arrow{ne,b}{0_{0,Y}}
\end{diagram}
$$
Certainly, such a morphism is also unique, and it has a special notation
$$
\ph=0_{X,Y}.
$$
In each category with a zero object the notions of {\it kernel} and {\it cokernel} are naturally defined \cite[Chapter VIII, $\S$ 1]{MacLane}.

\bex A key example for us is the category of augmented algebras. Recall that an {\it augmentation} on a (unital and associative) algebra $A$ over $\C$ is an arbitrary homomorphism  $\e:A\to\C$ of unital algebras over $\C$. It is easy to see, that the choice of the augmentation on $A$ is equivalent to the choice of a two-sided ideal $I_A$ in $A$ such that $A$ is a direct sum of the spaces over $\C$
$$
A=I_A\oplus \C\cdot 1_A,
$$
where $1_A$ is the unit of $A$.

An {\it augmented algebra} is a pair $(A,\e)$, where $A$ is a unital algebra over $\C$, and $\e:A\to\C$ an augmentation. The class of all augmented algebras forms a category $\Aug$, where morphisms $\ph:(A,\e_A)\to (B,\e_B)$ are arbitrary homomorphisms $\ph:A\to B$ of unital algebras over $\C$, which preserve augmentation in the following sense:
$$
\e_B\circ\ph=\e_A.
$$
We note that
\bit{\it
\item[1)] the algebra $\C$ with the identity mapping $\id_{\C}:\C\to\C$ as an augmentation is a zero object in the category $\Aug$.

\item[2)] the category $\Aug$ of augmented algebras over $\C$ has kernels and cokernels.
}\eit

\eex

\subsection{Kernel and cokernel as functors.}

Let $\tt K$ be an arbitrary category with a zero object. For each morphism $\alpha:A\to A'$ in this category we denote by $\Dom(\alpha)$ and $\Ran(\alpha)$ respectively its domain and range (codomain):
$$
\Dom(\alpha)=A,\qquad \Ran(\alpha)=A'.
$$
Thus each morphism $\alpha$ goes from $\Dom(\alpha)$ to $\Ran(\alpha)$:
$$
\alpha: \Dom(\alpha)\to \Ran(\alpha).
$$
Let us form a new category $\Arr({\tt K})$ from $\tt K$ by the following rules:

\bit{

\item[---] the objects of $\Arr({\tt K})$ are morphisms of the category ${\tt K}$:
$$
\Ob(\Arr({\tt K}))=\Mor({\tt K}),
$$

\item[---] a morphism in the category $\Arr({\tt K})$ between objects  $\alpha,\beta\in\Ob(\Arr({\tt K}))=\Mor(\tt K)$ is an arbitrary pair $(\ph,\psi)$ of morphisms $\ph,\psi\in\Ob(\Arr({\tt K}))=\Mor(\tt K)$ such that the following diagram is commutative:
 \beq\label{morfism-v-K^*}
 \xymatrix @R=2.pc @C=4.pc
 {
 \Dom(\alpha)\ar[d]^{\alpha}  \ar@{-->}[r]^{\ph} & \Dom(\beta) \ar[d]^{\beta}
 \\
 \Ran(\alpha) \ar@{-->}[r]^{\psi} & \Ran(\beta)
 }
 \eeq

\item[---] a composition of morphisms $(\ph,\psi)$ and $(\chi,\omega)$ in $\Arr({\tt K})$ is the pair $(\chi\circ\ph,\omega\circ\psi)$ (under the assumption that the compositions  $\chi\circ\ph$ and $\omega\circ\psi$ are defined in ${\tt K}$); this is illustrated by the diagram
 \beq\label{kompoz-v-K^*}
 \xymatrix @R=2.pc @C=4.pc 
 {
 \Dom(\alpha)\ar[d]^{\alpha}  \ar@{-->}[r]^{\ph} & \Dom(\beta) \ar[d]^{\beta}\ar@{-->}[r]^{\chi} & \Dom(\gamma) \ar[d]^{\gamma}
 \\
 \Ran(\alpha) \ar@{-->}[r]^{\psi} & \Ran(\beta)\ar@{-->}[r]^{\omega} & \Ran(\gamma)
 }
 \eeq

}\eit

The category $\Arr({\tt K})$ is called the {\it category of arrows} of the category ${\tt K}$ \cite[Chapter II, $\S$ 4]{MacLane}. Apparently, it is convenient to write the morphisms in the category $\Arr({\tt K})$, as special fractions: if the morphism $(\ph,\psi):\alpha\to\beta$ in \eqref{morfism-v-K^*} is represented by the symbol
$$
\frac{\beta,\ph}{\psi,\alpha},
$$
then the law of composition in $\Arr({\tt K})$ will be presented by the formula
\beq
\frac{\gamma,\chi}{\omega,\beta}\circ\frac{\beta,\ph}{\psi,\alpha}=\frac{\gamma,\chi\circ\ph}{\omega\circ\psi,\alpha}.
\eeq

The following proposition is obvious.

\btm\label{TH:Ker-frac-beta,ph-ph',alpha}
If the category $\tt K$ has kernels, then
\bit{
\item[(i)] in the category $\Arr({\tt K})$ each morphism $\frac{\beta,\ph}{\psi,\alpha}$ generates a unique morphism $\Ker\frac{\beta,\ph}{\psi,\alpha}: \Ker\alpha\to\Ker\beta$ such that the following diagram in $\tt K$ is commutative:
 \beq\label{Ker-frac-beta,ph-ph',alpha}
 \xymatrix @R=2.pc @C=4.pc
 {
 \Ker\alpha\ar[d]_{\ker\alpha}\ar@{-->}[r]^{\Ker\frac{\beta,\ph}{\psi,\alpha}} & \Ker\beta\ar[d]_{\ker\beta} \\
 \Dom(\alpha)\ar[d]^{\alpha}  \ar[r]^{\ph} & \Dom(\beta) \ar[d]^{\beta}
 \\
 \Ran(\alpha) \ar[r]^{\psi} & \Ran(\beta)
 }
 \eeq

 \item[(ii)] the mapping
 $$
\frac{\beta,\ph}{\psi,\alpha}\mapsto \Ker\frac{\beta,\ph}{\psi,\alpha}
 $$
is a covariant functor from $\Arr({\tt K})$ into ${\tt K}$;

 \item[(iii)] the mapping
 $$
\frac{\beta,\ph}{\psi,\alpha}\mapsto \frac{\ker\beta,\Ker\frac{\beta,\ph}{\psi,\alpha}}{\ph,\ker\alpha}
 $$
is a covariant functor from $\Arr({\tt K})$ into $\Arr({\tt K})$.

 }\eit
\etm

\bprop\label{TH:Ker-frac-beta,ph-ph',alpha-in-Mono}
If in the diagram \eqref{Ker-frac-beta,ph-ph',alpha} the morphism $\ph$ is a monomorphism, then the morphism $\Ker\frac{\beta,\ph}{\psi,\alpha}$ is a monomorphism as well;
\eprop
\bpr
The morphism
$$
\ker\beta\circ\Ker\frac{\beta,\ph}{\psi,\alpha}=\ph\circ\ker\alpha
$$
is a monomorphism as a composition of two monomorphisms, $\ker\alpha$ and $\ph$, hence the inner morphism in the left composition, $\Ker\frac{\beta,\ph}{\psi,\alpha}$, is a monomorphism as well.
\epr

For the cokernels the proposition dual to Theorem \ref{TH:Ker-frac-beta,ph-ph',alpha} is the following:

\btm\label{TH:Coker-frac-beta,ph-ph',alpha}
If the category $\tt K$ has cokernels, then
\bit{
\item[(i)] in the category $\Arr({\tt K})$ each morphism $\frac{\beta,\ph}{\psi,\alpha}$ generates a unique morphism $\Coker\frac{\beta,\ph}{\psi,\alpha}: \Coker\alpha\to\Coker\beta$ such that the following diagram in $\tt K$ is commutative:
 \beq\label{Coker-frac-beta,ph-ph',alpha}
 \xymatrix @R=2.pc @C=4.pc
 {
 \Dom(\alpha)\ar[d]_{\alpha}  \ar[r]^{\ph} & \Dom(\beta) \ar[d]^{\beta}
 \\
 \Ran(\alpha) \ar[r]^{\psi}\ar[d]_{\coker\alpha} & \Ran(\beta) \ar[d]^{\coker\beta} \\
  \Coker\alpha\ar@{-->}[r]^{\Coker\frac{\beta,\ph}{\psi,\alpha}} & \Coker\beta
 }
 \eeq

 \item[(ii)] the mapping
 $$
\frac{\beta,\ph}{\psi,\alpha}\mapsto \Coker\frac{\beta,\ph}{\psi,\alpha}
 $$
is a covariant functor from $\Arr({\tt K})$ into ${\tt K}$;

 \item[(iii)] the mapping
 $$
\frac{\beta,\ph}{\psi,\alpha}\mapsto \frac{\coker\beta,\psi}{\Coker\frac{\beta,\ph}{\psi,\alpha},\coker\alpha}
 $$
is a covariant functor from $\Arr({\tt K})$ into $\Arr({\tt K})$.

 }\eit
\etm

And the proposition dual to Proposition \ref{TH:Ker-frac-beta,ph-ph',alpha-in-Mono} is

\bprop\label{TH:Coker-frac-beta,ph-ph',alpha-in-Epi}
If in Diagram \eqref{Coker-frac-beta,ph-ph',alpha} the morphism $\psi$ is an epimorphism, then the morphism $\Coker\frac{\beta,\ph}{\psi,\alpha}$ is an epimorphism as well.
\eprop

\section{Stereotype spaces}

We shall need some definitions and facts from the theory of stereotype spaces (see details in \cite{Akbarov,Akbarov-stein-groups,Akbarov-env,Akbarov-C^infty-1}). We use the abbreviation  ``LCS'' for (Hausdorff) locally convex spaces over the field of complex numbers $\C$ \cite{Schaefer,Jarchow}.

A set $S$ in an LCS $X$ is said to be {\it totally bounded}, if for each neighbourhood of zero $U$ in $X$ there is a finite set $A$ such that the shifts of $U$ by elements of $A$ cover $S$:
$$
S\subseteq U+A.
$$
This is equivalent to the total boundedness of $S$ in the sense of the uniform structure induced from $X$ \cite{Engelking} (i.e.  {\it $A$ can be chosen from $S$}).

An LCS $X$ is said to be {\it pseudocomplete}, if each totally bounded Cauchy net in $X$ converges. This is equivalent to the condition that each closed totally bounded set $S$ in $X$ is compact. This notion is related to the usual completeness and quasicompleteness\footnote{An LCS $X$ is said to be {\it quasicomplete}, if each bounded Cauchy net in $X$ converges.} by the chain of implications
\beq\label{implikatsii-dlya-psevdopolnoty}
\text{$X$ is complete $\Longrightarrow$ $X$ is quasicomplete $\Longrightarrow$ $X$ is pseudocomplete.}
\eeq
In the metrizable case these conditions are equivalent.

Like completeness, pseudocompleteness has the property that each LCS $X$ has an ``outside-nearest'' pseudocomplete space. Formally this construction is described in the following

\btm\label{TH:DEF:triangledown} \cite[Theorem 1.5]{Akbarov} For each LCS $X$ there exists a pseudocomplete LCS $X^\triangledown$ and a linear continuous mapping ${\triangledown}_X :X \to X^{\triangledown}$ such that for each pseudocomplete LCS $Y$ and for each linear continuous mapping $\varphi:X\to Y$ there is a unique linear continuous mapping $\varphi^\triangledown:X^\triangledown\to Y$ such that the following diagram is commutative:
\beq\label{DEF:triangledown}
 \xymatrix @R=2.pc @C=1.pc %
 {
 X\ar[rr]^{\triangledown_X}\ar[dr]_{\varphi}& & X^\triangledown\ar@{-->}[dl]^{\varphi^\triangledown} \\
 & Y &
 }
\eeq
\etm

The space $X^\triangledown$ is called a {\it pseudocompletion} of the space $X$. The fact that the arrow $\varphi^\triangledown$ in \eqref{DEF:triangledown} is unique implies that $X^\triangledown$ is defined uniquely up to an isomorphism. On the other hand, Theorem \ref{TH:DEF:triangledown} implies that for each linear continuous mapping of locally convex spaces $\ph :X\to Y$ there is a unique linear continuous mapping $\ph^{\triangledown} :X^{\triangledown} \to Y^{\triangledown}$ such that the following diagram is commutative:
\beq
\begin{diagram}
\node{X} \arrow{e,t}{\triangledown_X} \arrow{s,l}{\varphi}
\node{X^\triangledown} \arrow{s,r,--}{\varphi^\triangledown}
\\
\node{Y} \arrow{e,t}{\triangledown_Y} \node{Y^\triangledown}
\end{diagram}.
\label{eqI.12}
\eeq
Moreover, the correspondence $\varphi\mapsto\varphi^\triangledown$ can be defined as a functor (i.e. as a mapping $\varphi\mapsto\varphi^\triangledown$ of the class of morphisms of the category of locally convex spaces into itself \cite[Theorem 1.5]{Akbarov}).

The pseudocompletion $X^{\triangledown}$ can be conceived as an {\it envelope (in the sense of \cite{Akbarov-env}) of the locally convex space $X$ in the class $\tt PC$ of all pseudocomplete locally convex spaces with respect to the same class $\tt PC$}:
\beq\label{X^triangledown=Env_PC^PC-X}
X^{\triangledown}=\Env_{\tt PC}^{\tt PC} X=\Env^{\tt PC} X
\eeq
(this fact follows from Theorem \ref{TH:DEF:triangledown} and the remark at page 50 of the work \cite{Akbarov-env}). This operation is similar to the usual completion in that it does not change the topology of the space $X$ (i.e. the mapping ${\triangledown}_X :X \to X^{\triangledown}$ is not only continuous, but is injective and open into its image), but only adds densely some elements to $X$ (i.e. the mapping ${\triangledown}_X$ embeds $X$ densely into $X^{\triangledown}$).

A set $D$ in a locally convex space $X$ is said to be {\it capacious}, if for each totally bounded set $S\subseteq X$ there is a finite set $A\subseteq X$ such that the shifts of $D$ by elements of $A$ cover $S$:
$$
S\subseteq D+A.
$$
It is useful to note that {\it if $D$ is convex, then $A$ can be chosen in $S$} (and this leads to an equivalent condition for $D$).

An LCS $X$ is said to be {\it pseudosaturated}, if each closed convex balanced set $D$ in $X$ is a neibourhood of zero in $X$. In the theory of topological vector spaces this property is connected with metrizability and barreledness by the following implications:
\beq\label{implikatsii-dlya-psevdonasyshennosti}
\text{$X$ is metrizable $\Longrightarrow$ $X$ is barreled $\Longrightarrow$ $X$ is pseudosaturated}
\eeq

It is remarkable that there is a standard construction, dual in a certain sense to the construction of pseudocompleteness, that allows to each LCS $X$ assign an ``inside-nearest'' pseudosaturated LCS $X^\vartriangle$:

\btm\label{TH:DEF:vartriangle} \cite[Theorem 1.16]{Akbarov} For each LCS $X$ there is a pseudosaturated LCS $X^\vartriangle$ and a linear continuous mapping ${\vartriangle}_X :X^{\vartriangle} \to X$ such that for each pseudosaturated LCS $Y$ and for each linear continuous mapping $\varphi:Y\to X$ there is a unique linear continuous mapping $\varphi^\vartriangle:Y\to X^\vartriangle$ such that the following diagram is commutative:
\beq\label{DEF:vartriangle}
 \xymatrix @R=2.pc @C=1.pc %
 {
 X& & X^\vartriangle\ar[ll]_{\vartriangle_X} \\
 & Y\ar[ul]^{\varphi} \ar@{-->}[ur]_{\varphi^\vartriangle}&
 }
\eeq
\etm

The space $X^\vartriangle$ is called a {\it pseudosaturation} of the space $X$. The fact that the arrow $\varphi^\vartriangle$ in \eqref{DEF:vartriangle} is unique implies that the LCS $X^\vartriangle$ is defined uniquely up to an isomorphism. On the other hand, Theorem \ref{TH:DEF:vartriangle} implies that for each linear continuous mapping of locally convex spaces $\ph:Y\to X$ there is a unique linear continuous mapping $\ph^\vartriangle :Y^\vartriangle\to X^\vartriangle$ such that the following diagram is commutative:
\beq
\begin{diagram}
\node{X} \node{X^\vartriangle} \arrow{w,t}{\vartriangle_X}
\\
\node{Y} \arrow{n,l}{\varphi} \node{Y^\vartriangle} \arrow{w,t}{\vartriangle_Y}
\arrow{n,r,--}{\varphi^\vartriangle}
\end{diagram}.
\label{eqI.25}
\eeq
Moreover, the correspondence $\varphi\mapsto\varphi^\vartriangle$ can be defined as a functor (i.e. as a mapping $\varphi\mapsto\varphi^\triangledown$ of the class of morphisms of the category of locally convex spaces into itself \cite[Theorem 1.16]{Akbarov}).

The pseudosaturation $X^{\vartriangle}$ can be conceived as a {\it refinement (in the sense of \cite{Akbarov-env}) of the locally convex space $X$ in the class $\tt PS$ of all pseudosaturated locally convex spaces by means of the same class $\tt PS$}:
\beq\label{X^vartriangle=Ref_PS^PS-X}
X^{\vartriangle}=\Rf_{\tt PS}^{\tt PS} X=\Rf^{\tt PS} X
\eeq
(this fact follows from Theorem \ref{TH:DEF:vartriangle} and the remark at pages 56-57 in \cite{Akbarov-env}). This operation does not change the composition of elements of $X$ (i.e. the mapping ${\vartriangle}_X :X^{\vartriangle}\to X$ is bijective), but strengthen the topology of $X$ in some special way such that the system of totally bounded sets in $X$ is not changed, and the topology on each totally bounded set $S\subseteq X$ is not changed as well.

A {\it dual space} to an LCS $X$ is defined as the space $X^\star$\label{DEF:X^star} of all linear continuous functionals $f:X\to \C$ endowed with the {\it topology of uniform convergence on totally bounded sets} in $X$. The {\it second dual space} $X^{\star\star}$ is the dual to $X^\star$ in the same sense:
$$
  X^{\star\star}=(X^\star)^\star.
$$
The formula
\beq\label{DEF:i_X(x)(f)=f(x)}
     i_X(x)(f)=f(x),\qquad x\in X,
\eeq
defines a natural mapping $i_X:X\to X^{\star\star}$.

A locally convex space  $X$ is said to be {\it stereotype}, if the mapping $i_X:X\to X^{\star\star}$ is an isomorphism of locally convex spaces (i.e. a homeomorphism between $X$ and $X^{\star\star}$).

\btm \cite[Theorem 4.1]{Akbarov}
A locally convex space $X$ is stereotype if and only if it is pseudocomplete and pseudosaturated.
\etm

The class $\tt Ste$ of stereotype spaces is extremely wide, it contains, in particular, all pseudocomplete barreled spaces (and hence all Banach spaces and all Fr\'echet spaces).

\section{Kernels and cokernels in the category of augmented involutive stereotype algebras.}

A stereotype space $A$ (over $\C$) is called a {\it stereotype algebra}, if $A$ is endowed with the structure of a unital associative algebra (over $\C$), and the multiplication operation is a continuous map in the following sense: for each compact set $K$ in $A$ and for any neighborhood of zero $U$ in $A$ there is a neighborhood of zero $V$ in $A$ such that
$$
K\cdot V\subseteq U\quad \& \quad V\cdot K\subseteq U.
$$
An {\it involution} on a stereotype algebra $A$ is an arbitrary continuous mapping $\bullet:A\to A$ with the properties
$$
(x+y)^\bullet=x^\bullet+y^\bullet, \quad
(\lambda\cdot x)^\bullet=\overline{\lambda}\cdot x^\bullet, \quad
(x\cdot y)^\bullet=y^\bullet\cdot x^\bullet,\quad
(x^\bullet)^\bullet=x,\quad x,y\in A,\ \lambda\in\C.
$$
A pair $(A,\bullet)$, where $A$ is a stereotype algebra, and $\bullet$ an involution on $A$, is called an {\it involutive stereotype algebra}.

The class of all involutive stereotype algebras is denoted by $\InvSteAlg$. It forms a category where morphisms are the continuous linear multiplicative and unit- and involution-preserving maps $\ph :A\to B$.

An {\it augmentation} on an involutive stereotype algebra $A$ is defined as a morphism $\e:A\to\C$ in $\InvSteAlg$. An {\it augmented involutive stereotype algebra} is a pair $(A,\e)$, where $A$ is an involutive stereotype algebra over $\C$, and $\e:A\to\C$ an augmentation on it. Defining an augmentation on a stereotype algebra $A$ is equivalent to defining a two-sided closed and invariant with respect to the involution ideal $I_A$ in $A$ such that $A$ is a direct sum of the vector spaces over $\C$
$$
A=I_A\oplus \C\cdot 1_A,
$$
where $1_A$ is the unit in $A$.

The class of all augmented involutive stereotype algebras will be denoted by $\AugInvSteAlg$. It forms a category where morphisms $\ph :(A,\e_A)\to (B,\e_B)$ are morphisms of stereotype algebras $\ph :A\to B$ which preserve the involution and the augmentation: $\e_A=\e_B\circ\ph$.

As in the category $\Aug$, the algebra $\C$ with the identity map $\id_{\C}:\C\to\C$ as an augmentation is the zero object in $\AugInvSteAlg$. As a corollary, a zero morphism in  $\AugInvSteAlg$ is an arbitrary morphism $\ph:A\to B$ which can be decomposed as $\iota_B\circ\e_A$, where $\e_A:A\to\C$ is the augmentation on $A$, and $\iota_B:\C\to B$ is the  (unique) morphism of $\C$ into $B$. Thus,
\beq\label{0=iota_B-circ-e_A}
0_{A,B}=\iota_B\circ\e_A.
\eeq

\bprop\label{PROP:0:A->B}
A morphism $\ph:A\to B$ of augmented involutive stereotype algebras is a zero morphism if and only if its set of values is contained in the subalgebra generated by the unit $1_B$ of the algebra $B$:
\beq\label{0:A->B}
\ph=0_{A,B}\quad\Longleftrightarrow\quad \ph(A)\subseteq\C\cdot 1_B.
\eeq
\eprop

\bex\label{EX:functional-algebras} {\sf Functional algebras on groups \cite{Akbarov-group-algebras,Akbarov}.} A typical example of an augmented involutive stereotype algebra is the algebra ${\mathcal C}(G)$ of continuous functions on a locally compact group $G$ with the pointwise algebraic operations and the topology of uniform convergence on compact sets in $G$. The fact that it is a stereotype algebra is proved in \cite[Example 10.3]{Akbarov}. The involution and the augmentation on ${\mathcal C}(G)$ are defined by the formulas
\beq\label{u^bullet(t)=overline{u(t)}}
u^\bullet(t)=\overline{u(t)},\qquad \e(u)=u(1_G),\quad u\in {\mathcal C}(G).
\eeq
This example can be complemented by another classical functional algebra on groups, the algebra ${\mathcal E}(G)$ of smooth functions on a Lie group $G$ (with the usual topology of uniform convergence on compact sets with respect to each partial derivative \cite[Example 10.4]{Akbarov}). The involution and the augmentation on ${\mathcal E}(G)$ are defined by the same formulas \eqref{u^bullet(t)=overline{u(t)}}.
\eex

\bex\label{EX:group-algebras} {\sf Group algebras \cite{Akbarov-group-algebras,Akbarov}.}
When we take dual spaces we obtain two other examples, the algebras with respect to convolution:
\bit{

\item[---] the algebra ${\mathcal C}^\star(G)$ of measures with compact support on a locally compact group $G$ \cite[Example 10.7]{Akbarov}, and

\item[---] the algebra ${\mathcal E}^\star(G)$ of distributions with compact support on a Lie group $G$ (with the usual topology of uniform convergence on compact sets in ${\mathcal E}(G)$ \cite[Example 10.8]{Akbarov}).

}\eit
The involution in these algebras is defined by the formula
 \beq\label{involution-in-C*(G)}
\alpha^\bullet(u)=\overline{\alpha(\widetilde{u}^\bullet)}, \qquad \alpha\in {\mathcal C}^\star(G)\quad ({\mathcal E}^\star(G)),
 \eeq
where $u^\bullet$ is the involution defined by the identity \eqref{u^bullet(t)=overline{u(t)}}, and
$$
\widetilde{u}(t)=u(t^{-1})
$$
is the antipode of the function. The augmentation is
 \beq\label{augmentation-in-C*(G)}
\e(\alpha)=\alpha(1),\qquad \alpha\in {\mathcal C}^\star(G)\quad ({\mathcal E}^\star(G))
 \eeq
($1$ --- is the function on $G$ identically equal to 1).
\eex

\subsection{Existence of kernels and cokernels in $\AugInvSteAlg$.}

Further we shall need the following constructions.
\bit{
\item Let $B$ be a stereotype algebra and $A$ a closed subalgebra in $B$ (i.e. $A$ is a unital subalgebra in $B$ in the purely algebraic sense and at the same time a closed subspace in the locally convex space $B$). Let us endow $A$ with the topology induced by $B$. Then the pseudosaturation $A^\vartriangle$ of the space $A$ is a stereotype subalgebra, called a {\it closed immediate subalgebra in the stereotype algebra $B$, generated by the subalgebra $A$} \cite[Theorem 5.14]{Akbarov-env}.

\item Again, let $B$ be a stereotype algebra, and suppose $I$ is a closed two-sided ideal in $B$ (i.e. $I$ is a two-sided ideal in $B$ in the purely algebraic sense and at the same time a closed subspace in the locally convex space $B$). Consider the quotient space $B/I$. It is an algebra in the purely algebraic sense, as the quotient algebra of the algebra $B$ by the ideal $I$, and at the same time $B/I$ is a locally convex space with the usual quotient topology inherited from $B$. The pseudocompletion $(B/I)^\triangledown$ of the space $B/I$ is a stereotype algebra in $B$, called the {\it open immediate quotient algebra of the stereotype algebra $B$ by the ideal $I$} \cite[Theorem 5.17]{Akbarov-env}.
}\eit

\btm\label{PROP:AugSteAlg-imeet-ker}
For a morphism $\ph:(A,\e_A)\to (B,\e_B)$ in the category $\AugInvSteAlg$ of augmented involutive stereotype algebras over $\C$
\bit{

\item[---] the kernel is the closed immediate subalgebra in $A$, generated by the inverse image of the subalgebra $\C\cdot 1_B$ in $B$ under the mapping $\ph$:
\beq\label{ker-v-AugSteAlg}
\Ker\ph=\big(\ph^{-1}(\C\cdot 1_B)\big)^\vartriangle,
\eeq

\item[---] the cokernel is the open immediate quotient algebra of the algebra $B$ by the closed two-sided ideal $I$ in $B$, generated by the image $\ph(\e_A^{-1}(0))$ of the ideal $\e_A^{-1}(0)$ in $A$:
\beq\label{coker-v-AugSteAlg}
\Coker\ph=\big(B/I\big)^\triangledown.
\eeq
}\eit
As a corollary, the category $\AugInvSteAlg$ has kernels and cokernels.
\etm
\bpr
1. Let us prove \eqref{ker-v-AugSteAlg}. Set $K=\big(\ph^{-1}(\C\cdot 1_B)\big)^\vartriangle$ (we endow $K$ with the structure of an immediate subspace in $A$). Let us show that $\varkappa:K\to A$ is a kernel of the morphism $\ph:A\to B$.

a) First, from Proposition \ref{PROP:0:A->B} it follows that $\ph\circ\varkappa:K\to B$ is zero, since
$$
\varkappa(K)\subseteq \C\cdot 1_B.
$$

b) Second, let $\varkappa':(K',\e_{K'})\to (A,\e_A)$ be another morphism that gives zero in composition with $\ph$:
$$
\ph\circ\varkappa'=0=\iota_B\circ\e_{K'}.
$$
Then
$$
\ph(\varkappa'(K'))\subseteq\Im\iota_B=\C\cdot 1_B\quad\Rightarrow\quad
\varkappa'(K')\subseteq\ph^{-1}(\C\cdot 1_B)=K,
$$
hence $\varkappa'$ is lifted to a morphism of stereotype spaces $\delta:K'\to K$. It is a morphism of augmented involutive stereotype algebras, since it preserves multiplication, unit, involution and augmentation.

2. Let us prove \eqref{coker-v-AugSteAlg}. We endow the ideals
$$
I_A=\e_A^{-1}(0),\qquad I_B=\e_B^{-1}(0)
$$
with the structure of immediate subspaces in $A$ and $B$ (in this case with the topology which is a pseudosaturation of the topology induced from $A$ and $B$, see details in  \cite[Proposition 4.69]{Akbarov-env}). Let $I$ be a closed immediate two-sided ideal in $B$ generated by the set $\ph(I_A)$ (i.e. $I$, as a set, coincides with the closure in $B$ of the two-sided ideal generated by the set $\ph(I_A)$, and is endowed with the topology which is a pseudosaturation of the topology induced from $B$). The equality $\e_B\circ\ph=\e_A$ implies
$$
\ph(I_A)\subseteq I_B.
$$
Therefore,
$$
\ph(I_A)\subseteq I\subseteq I_B.
$$
Put $C=(B/I)^\triangledown$, and denote by $\gamma:B\to (B/I)^\triangledown=C$ the quotient map. The condition $I\subseteq I_B=\e_B^{-1}(0)$ implies that the functional  $\e_B$ can be extended to some functional $\e_C$ on the quotient space $(B/I)^\triangledown=C$:
$$
 \xymatrix @R=2.pc @C=1.pc %
 {
 B\ar[rr]^{\gamma}\ar[dr]_{\e_B}& & C\ar[dl]^{\e_C} \\
 & \C &
 }
$$
Hence $\gamma$ is a morphism of augmented stereotype algebras $(B,\e_B)\to(C,\e_C)$. Let us show that this is a cokernel of the morphism $\ph$.

a) First, $\gamma\circ\ph$ is zero. To show this take an $a\in A$ and put  $x=a-\e_A(a)\cdot 1_A$. Then
\begin{eqnarray*}
a\in A\quad\Longrightarrow\quad a=\e_A(a)\cdot 1_A+\underset{\scriptsize\begin{matrix}\text{\rotatebox{90}{$\owns$}}\\ I_A\end{matrix}}{x}\quad\Longrightarrow\quad
\ph(a)=\e_A(a)\cdot 1_B+\underset{\scriptsize\begin{matrix}\text{\rotatebox{90}{$\owns$}}\\ I\end{matrix}}{\ph(x)}\quad\Longrightarrow \\ \Longrightarrow\quad
\gamma(\ph(a))=\e_A(a)\cdot 1_C+\underbrace{\gamma(\ph(x))}_{\scriptsize\begin{matrix}\text{\rotatebox{90}{$=$}}\\ 0\end{matrix}}=\e_A(a)\cdot 1_C\in\C\cdot 1_C.
\end{eqnarray*}
So $(\gamma\circ\ph)(A)\subseteq \C\cdot 1_C$, and by Proposition \ref{PROP:0:A->B} this means that $\gamma\circ\ph=0$.

b) Now let $\gamma':(B,\e_B)\to(C',\e_{C'})$ be another morphism that gives zero in composition with $\ph$:
\beq\label{delta-circ-ph=iota_D-circ-e_A}
\gamma'\circ\ph=\iota_{C'}\circ\e_A.
\eeq
To verify that it can be factored through $\gamma$, i.e. for some morphism $\iota$ the diagram
$$
 \xymatrix @R=2.pc @C=4.pc
 {
 A\ar[r]^{\ph}\ar[dr]_{0_{A,C'}} & B\ar[d]_{\gamma'} \ar[r]^{\gamma}& C \ar@{-->}[dl]^{\iota} \\
 & C' &
 }
$$
is commutative, we have to verify that
$$
I\subseteq(\gamma')^{-1}(0).
$$
In other words,
$$
\gamma'(I)=0.
$$
Since $I$ is a closed two-sided ideal generated by the set $\ph(I_A)$, it is sufficient for us to prove that
$$
\gamma'(\ph(I_A))=0.
$$
This follows from \eqref{delta-circ-ph=iota_D-circ-e_A}:
$$
x\in I_A\quad\Longrightarrow\quad
\gamma'(\ph(x))=\eqref{delta-circ-ph=iota_D-circ-e_A}=\iota_{C'}(\e_A(\underset{\scriptsize\begin{matrix}\text{\rotatebox{90}{$\owns$}}\\ I_A\end{matrix}}{x}))=\iota_{C'}(0)=0.
$$
\epr

\subsection{Preservation of the cokernel by stereotype group algebras.}

Recall the group algebra ${\mathcal C}^\star(G)$ of measures from Example  \ref{EX:group-algebras}.

\btm\label{TH:coker-groups=>coker-group-algebras}
Suppose that in a chain of locally compact groups
$$
 \xymatrix  
{
H\ar[r]^{\lambda} & G\ar[r]^{\varkappa} &F
}
$$
the second homomorphism is a cokernel of the first one in the category of locally compact groups\footnote{In other words, $\varkappa$ is a quotient map of the group $G$ by the closure $\overline{\lambda(H)}$ of the subgroup $\lambda(H)$.}:
$$
\coker\lambda=\varkappa.
$$
Then in the corresponding chain of stereotype group algebras
$$
 \xymatrix 
{
{\mathcal C}^\star(H)\ar[r]^{{\mathcal C}^\star(\lambda)} & {\mathcal C}^\star(G)\ar[r]^{{\mathcal C}^\star(\varkappa)} &
{\mathcal C}^\star(F)
}
$$
the second morphism is a cokernel of the first one in the category $\AugInvSteAlg$ of augmented involutive stereotype algebras:
\beq\label{coker-C^star(lambda)=C^star(varkappa)}
\coker{\mathcal C}^\star(\lambda)={\mathcal C}^\star(\varkappa).
\eeq
\etm
\bpr
Let us introduce a simpler notations:
$$
\lambda'={\mathcal C}^\star(\lambda),\quad \varkappa'={\mathcal C}^\star(\varkappa).
$$
And let $\e=\e_{{\mathcal C}^\star(H)}$ be the augmentation on ${\mathcal C}^\star(H)$, and $\iota=\iota_{{\mathcal C}^\star(F)}$ the embedding of $\C$ into ${\mathcal C}^\star(F)$.

1. First let us note that $\varkappa'\circ\lambda'$ is a zero morphism:
\beq\label{varkappa'-circ-lambda'=iota-circ-e}
\varkappa'\circ\lambda'=\iota\circ\e.
\eeq
$$
 \xymatrix @R=2.pc @C=4.pc
{
H\ar[r]^{\lambda}\ar[d]_{\delta_H} & G\ar[r]^{\varkappa}\ar[d]_{\delta_G} &F\ar[d]_{\delta_F} \\
{\mathcal C}^\star(H)\ar[r]^{\lambda'}
\ar@{-->}@/_2ex/[rd]^{\e} & {\mathcal C}^\star(G)\ar[r]^{\varkappa'} &
{\mathcal C}^\star(F)\\
& \C \ar@{-->}@/_2ex/[ru]^{\iota} &
}
$$
It is convenient to see this on elements $\delta^t$, $t\in H$:
$$
(\varkappa'\circ\lambda')(\delta^t)=\varkappa'(\lambda'(\delta^t))=\varkappa'(\delta^{\lambda(t)})=\delta^{\varkappa(\lambda(t))}=
\delta^{1_F}=1_{{\mathcal C}^\star(F)}=\iota(1)=\iota(\e(\delta^t))=(\iota\circ\e)(\delta^t).
$$
Since elements $\delta^t$, $t\in H$, are full in ${\mathcal C}^\star(H)$ \cite[Lemma 8.2]{Akbarov}, this proves \eqref{varkappa'-circ-lambda'=iota-circ-e}.

2. Now let $\gamma:{\mathcal C}^\star(G)\to A$ be a morphism of augmented stereotype algebras that gives zero in composition with $\lambda'$:
\beq\label{gamma-circ-lambda'=iota-circ-e}
\gamma\circ\lambda'=\iota\circ\e.
\eeq
Consider the diagram:
$$
 \xymatrix @R=2.pc @C=4.pc
{
H\ar[r]^{\lambda}\ar[d]_{\delta_H} & G\ar[rr]^{\varkappa}\ar[d]_{\delta_G} & & F\ar[dl]_{\delta_F}\ar@{-->}@/^6ex/[ddll]^{\ph} \\
{\mathcal C}^\star(H)\ar[r]^{\lambda'}\ar[d]_{\e}
 & {\mathcal C}^\star(G)\ar[r]^{\varkappa'}\ar[d]_{\gamma} &
{\mathcal C}^\star(F)\ar@{-->}[dl]_{\gamma'}&  \\
\C\ar[r]^{\iota} & A & &
}
$$
For each $t\in H$ we have
$$
\gamma(\delta_G^{\lambda(t)})=
\gamma\big(\lambda'(\delta_H^t)\big)=\eqref{gamma-circ-lambda'=iota-circ-e}=
\iota\big(\e(\delta_H^t)\big)=\iota(1)=1.
$$
Thus, the homomorphism $G\to A$ is constant on the subgroup $\lambda(H)\subseteq G$. This means that it is induced by a homomorphism $\ph:F=\coker\lambda=G/\overline{\lambda(H)}\to A$. This homomorphism $\ph$ in its turn can be extended to a morphism of the group algebra $\gamma':{\mathcal C}^\star(F)\to A$ \cite[Theorem 10.12]{Akbarov}, which is the extension of the homomorphism $\gamma$, since on elements $G$ they coincide (we use here \cite[Lemma 8.2]{Akbarov}).
\epr

\subsection{Continuous envelope}

The notion of continuous envelope of an involutive stereotype algebra $A$ was introduced by the author in \cite{Akbarov-env} and was discussed in detail in \cite{Akbarov-C^infty-2}. This is an envelope of $A$ in the class $\DEpi$ of {\it dense epimorphisms}, i.e. the morphisms $\ph:A\to B$ with the property $\overline{\ph(A)}=B$, with the values in $C^*$-algebras:
$$
\Env_{\mathcal C}A=\Env_{\sc C^*}^{\DEpi}A.
$$
The detailed definition is the following.
 \bit{
\item First, a {\it continuous extension} of an involutive stereotype algebra $A$ is a dense epimorphism $\sigma:A\to A'$ (i.e. $\overline{\sigma(A)}=A'$) of involutive stereotype algebras such that for each $C^*$-algebra $B$ and for each involutive homomorphism $\ph:A\to B$ there is a (necessarily unique) homomorphism of involutive stereotype algebras $\ph':A'\to B$ such that the following diagram is commutative:
\beq\label{DEF:diagr-nepr-rasshirenie}
 \xymatrix @R=2pc @C=1.2pc
 {
  A\ar[rr]^{\sigma}\ar[dr]_{\ph} & & A'\ar@{-->}[dl]^{\ph'} \\
  & B &
 }
\eeq
\item  Second, a {\it continuous envelope} of an involutive stereotype algebra $A$ is a continuous extension $\env_{\mathcal C} A:A\to \Env_{\mathcal C} A$ such that for any continuous extension $\sigma:A\to A'$ there is a (necessarily unique) morphism of involutive stereotype algebras $\upsilon:A'\to \Env_{\mathcal C} A$ such that the following diagram is commutative:
$$
 \xymatrix @R=2pc @C=1.2pc
 {
  & A\ar[ld]_{\sigma}\ar[rd]^{\env_{\mathcal C} A} &   \\
  A'\ar@{-->}[rr]_{\upsilon} &  & \Env_{\mathcal C} A
 }
$$
 }\eit
If $\ph:A\to B$ is a morphism of involutive stereotype algebras, then there exists a unique morphism of involutive stereotype algebras $\Env_{\mathcal C} \ph:\Env_{\mathcal C} A\to \Env_{\mathcal C} B$ such that the following diagram is commutative \cite[(5.1.4)]{Akbarov-C^infty-2}:
    \beq\label{Env_C-ph}
     \xymatrix @R=2.pc @C=4.pc
 {
 A\ar[r]^{\env_{\mathcal C} A}\ar[d]_{\ph} &  \Env_{\mathcal C} A \ar@{-->}[d]^{\Env_{\mathcal C} \ph} \\
 B\ar[r]_{\env_{\mathcal C}B} & \Env_{\mathcal C}B
 }
    \eeq
The morphism $\Env_{\mathcal C} \ph$ is called the {\it continuous envelope of the morphism} $\ph$.

\btm\label{LM:A-augste=>Env_C-A-augste}
Let $(A,\e)$ be an augmented involutive stereotype algebra. Then
 \bit{

 \item[(i)] the continuous envelope $\Env_{\mathcal C} \e:\Env_{\mathcal C} A\to\Env_{\mathcal C} \C=\C$ of the augmentation $\e$ on $A$ is an augmentation on the continuous envelope $\Env_{\mathcal C} A$ of the algebra $A$;

 \item[(ii)] the envelope $\env_{\mathcal C} A:A\to\Env_{\mathcal C} A$ is a morphism of augmented involutive stereotype algebras.
    \beq\label{env_C-e}
     \xymatrix @R=2.pc @C=4.pc
 {
 A\ar[r]^{\env_{\mathcal C} A}\ar[d]_{\e} &  \Env_{\mathcal C} A \ar@{-->}[d]^{\Env_{\mathcal C} \e} \\
 \C\ar[r]_{\id_{\C}=\env_{\mathcal C} \C} & \C
 }
    \eeq
  }\eit
\etm
\bpr
This follows immediately from Diagram \eqref{env_C-e}, which is commutative due to \eqref{Env_C-ph} and \cite[(5.1.7)]{Akbarov-C^infty-2}.
\epr

\btm\label{TH:Env_C-functor-v-AugSteAlg}
If $\ph:(A,\e_A)\to(B,\e_B)$ is a morphism of augmented involutive stereotype algebras, then its continuous envelope $\Env_{\mathcal C} \ph:(\Env_{\mathcal C} A,\Env_{\mathcal C} \e_A)\to(\Env_{\mathcal C} B,\Env_{\mathcal C} \e_B)$ is a morphism of augmented involutive stereotype algebras as well due to the diagram
    \beq\label{Env_C-functor-v-AugSteAlg}
     \xymatrix @R=2.pc @C=4.pc
 {
 A\ar[rr]^{\env_{\mathcal C} A}\ar[dd]_{\ph}\ar[dr]_{\e_A} & &  \Env_{\mathcal C} A \ar@{-->}[dd]^{\Env_{\mathcal C} \ph}\ar@{-->}[dl]^{\ \Env_{\mathcal C} \e_A} \\
 & \C & \\
 B\ar[rr]_{\env_{\mathcal C} B}\ar[ur]^{\e_B} & & \Env_{\mathcal C} B\ar@{-->}[ul]_{\ \Env_{\mathcal C} \e_B}
 }
    \eeq
\etm
\bpr
Here the perimeter and all the inner triangles are commutative, excluding the right inner triangle (with the dashed arrows), and its commutativity we have to prove. But it is commutative since $\env_{\mathcal C} A$ is an epimorphism.
\epr

\bprop\label{LM:exten_C=>exten_C^Aug}
If $A$ is an involutive stereotype algebra with the augmentation $\e:A\to\C$, and $\sigma:A\to A'$ is its continuous extension, then the algebra $A'$ has a unique augmentation $\e':A'\to\C$ such that the morphism $\sigma$ becomes a morphism of augmented stereotype algebras, and moreover, an extension in the category $\AugInvSteAlg$ of augmented involutive stereotype algebras in the class $\DEpi$ of dense epimorphisms with respect to the class $\AugC*$ of the augmented $C^*$-algebras.
\eprop
\bpr
1. Consider the diagram:
    \beq\label{exten_C=>exten_C^Aug}
     \xymatrix 
 {
 A\ar[rr]^{\sigma}\ar[rd]_{\e} & & A' \ar@{-->}[dl]^{\e'} \\
 & \C &
 }
    \eeq
Since $\C$ is a $C^*$-algebra, the morphism $\e'$ exists and is uniquely defined. The commutativity of this diagram means that $\sigma$ is a morphism of augmented stereotype algebras $\sigma:(A,\e)\to(A',\e')$.

2. Let us show that this morphism $\sigma:(A,\e)\to(A',\e')$ is an extension in $\DEpi$ with respect to $\AugC*$. Let $\ph:(A,\e)\to(B,\delta)$ be a morphism into an augmented $C^*$-algebra. Consider a diagram in the category of stereotype algebras:
$$
     \xymatrix 
 {
 A\ar@/_5ex/[rdd]_{\ph}\ar[rr]^{\sigma}\ar[rd]_{\e} & & A' \ar[dl]^{\e'}\ar@{-->}@/^5ex/[ldd]^{\ph'} \\
 & \C & \\
 & B\ar[u]^{\delta} &
 }
$$
The dashed arrow, $\ph'$, exists, is unique and the perimeter becomes commutative, since  $B$ is a $C^*$-algebra, and $\sigma$ is a continuous envelope. At the same time the upper inner triangle is commutative, since this is just diagram \eqref{exten_C=>exten_C^Aug}, and the left inner triangle is commutative since $\ph:(A,\e)\to(B,\delta)$ is a morphism of augmented stereotype algebras. In addition $\sigma$ is an epimorphism, hence the right inner triangle is commutative as well:
$$
\delta\circ\ph'=\e'.
$$
This means that $\ph'$ is a morphism in the category $\AugInvSteAlg$, and since it is unique, $\sigma$ is an extension in $\AugInvSteAlg$ (in the class $\DEpi$ with respect to the class $\AugC*$ of augmented $C^*$-algebras).
\epr

\bcor
For each stereotype algebra $A$ with an augmentation $\e:A\to\C$ there is a unique morphism of augmented stereotype algebras $\upsilon:(A,\e)\to\Env_{\AugC*}^{\DEpi}(A,\e)$ such that the following diagram is commutative:
    \beq\label{DEF:Env_C-e}
     \xymatrix 
 {
 & (A,\e)\ar[dl]_{\env_{\mathcal C} A}\ar[dr]^{\env_{\AugC*}^{\DEpi}(A,\e)} & \\
 (\Env_{\mathcal C} A,\Env_{\mathcal C} \e) \ar@{-->}[rr]_{\upsilon}  & & \Env_{\AugC*}^{\DEpi}(A,\e)
 }
    \eeq
\ecor
\bpr
In the notations of Proposition \ref{LM:exten_C=>exten_C^Aug}, here $\sigma$ is not just an extension, but an envelope $\env_{\mathcal C} A$. In this case the augmentation $\e'$ on $A'=\Env_{\mathcal C} A$ is a morphism $\Env_{\mathcal C} \e$ from \eqref{env_C-e}. Since by Proposition \ref{LM:exten_C=>exten_C^Aug} $\sigma=\env_{\mathcal C} A$ is an extension in the category of augmented involutive stereotype algebras, there must be a unique morphism $\upsilon$ in \eqref{DEF:Env_C-e}.
\epr

\subsection{Continuous envelope of a group algebra ${\mathcal C}^\star(Z\cdot K)$}

The following fact was proved in \cite[Theorem 5.53]{Akbarov-env} (for the Kuznetsova envelopes in \cite[Theorem 2.11]{Kuznetsova}).

\btm\label{TH:Env_C-C^star(G)=C(widehat(G))}
The Fourier transform on an abelian locally compact group $Z$
\beq\label{Fourier:C^star(H)->C(widehat(H))}
{\mathcal F}_Z:{\mathcal C}^\star(Z)\to{\mathcal C}(\widehat{Z}),\quad
{\mathcal F}_Z(\alpha)(\chi)=\alpha(\chi),\quad  \alpha\in{\mathcal C}^\star(Z),\ \chi\in\widehat{Z},
\eeq
is a continuous envelope of the group algebra ${\mathcal C}^\star(Z)$. As a corollary,
\beq\label{Env_C-C^star(G)=C(widehat(G))}
\Env_{\mathcal C} {\mathcal C}^\star(Z)={\mathcal C}(\widehat{Z}).
\eeq
\etm

The following theorem is proved in \cite[Proposition 5.27]{Akbarov-C^infty-2}.

\btm\label{TH:env_C^star(Z-times-K)}
Suppose $Z$ is an abelian locally compact group, and $K$ a compact group. Then the formula
\beq\label{DEF:env_C^star(Z-times-K)}
(\Phi\delta^{(t,x)})_\sigma(\chi)=\chi(t)\cdot\sigma(x),\quad t\in Z,\ x\in K,\ \chi\in\widehat{Z}, \ \sigma\in \widehat{K},
\eeq
defines a mapping
\beq\label{env_C^star(Z-times-K)}
\Phi:{\mathcal C}^\star(Z\times K)\to\prod_{\sigma\in\widehat{K}}{\mathcal C}\big(\widehat{Z},{\mathcal B}(X_\sigma)\big),
\eeq
which is a continuous envelope of the group algebra ${\mathcal C}^\star(Z\times K)$. As a corollary,
\beq\label{env_C-C^star(Z-times-K)}
\Env_{\mathcal C} {\mathcal C}^\star(Z\times K)={\mathcal C}\Big(\widehat{Z},\prod_{\sigma\in\widehat{K}}{\mathcal B}(X_\sigma)\Big)=
\prod_{\sigma\in\widehat{K}}{\mathcal C}\big(\widehat{Z},{\mathcal B}(X_\sigma)\big).
\eeq
\etm

Let us call a locally compact group $G$ a {\it compact buildup of an abelian locally compact group}\label{DEF:komp-nadstr-abelevoi-gruppy}, if there exist closed subgroups $Z$ and $K$ in $G$ with the following properties:
\bit{

\item[1)] $Z$ is an abelian group,

\item[2)] $K$ is a compact group,

\item[3)] $Z$ and $K$ commute:
$$
\forall a\in Z,\quad \forall y\in K\qquad a\cdot y=y\cdot a,
$$

\item[4)] the product of $Z$ and $K$ is $G$:
$$
\forall x\in G\qquad \exists a\in Z\quad\exists y\in K\qquad x=a\cdot y.
$$
}\eit
If it is necessary to specify which groups in this constructions are used, then we say that $G$ is a {\it buildup of the abelian group $Z$ with the help of the compact group $K$}.

\blm\label{LM:Z-cdot-K=(Z-times-K)/C}
Let $G=Z\cdot K$ be a buildup of an abelian locally compact group $Z$ with the help of a compact group $K$. Consider the subgroup
$$
H=Z\cap K,
$$
its immersion into the cartesian product
$$
\iota:H\to Z\times K,\quad \iota(x)=(x,x^{-1}),\quad x\in H,
$$
and the image of this immersion
$$
\iota(H)\subseteq Z\times K.
$$
Then
$$
G\cong(Z\times K)/\iota(H).
$$
\elm

The homomorphism of the groups $\iota:H\to Z\times K$ generates a morphisms of group algebras ${\mathcal C}^\star(\iota):{\mathcal C}^\star(H)\to {\mathcal C}^\star(Z\times K)$, which in its turn generates a morphism of envelopes $\Env_{\mathcal C}{\mathcal C}^\star(\iota):\Env_{\mathcal C}{\mathcal C}^\star(H)\to \Env_{\mathcal C}{\mathcal C}^\star(Z\times K)$. As a corollary, we have the diagram
\beq\label{DEF:xi:Env-C*(H)->Env-C*(Z-times-K)}
 \xymatrix  @R=2.pc @C=6.pc
{
 H\ar[r]^{\iota}\ar[d]_{\delta_H} & Z\times K \ar[d]^{\delta_{Z\times K}}
 \\
{\mathcal C}^\star(H)\ar[r]^{{\mathcal C}^\star(\iota)}
 \ar[d]_{\env_{\mathcal C}{\mathcal C}^\star(H)} & {\mathcal C}^\star(Z\times K)
 \ar[d]^{\env_{\mathcal C}{\mathcal C}^\star(Z\times K)} \\
\Env_{\mathcal C}{\mathcal C}^\star(H)\ar[r]^{\Env_{\mathcal C}{\mathcal C}^\star(\iota)}\ar@{=}[d] & \Env_{\mathcal C}{\mathcal C}^\star(Z\times K)\ar@{=}[d]\\
{\mathcal C}(\widehat{H})\ar[r]^{\xi} & \prod_{\sigma\in\widehat{K}}{\mathcal C}\big(\widehat{Z},{\mathcal B}(X_\sigma)\big)
}
\eeq
where $\xi$ is the morphism, corresponding to $\Env_{\mathcal C}{\mathcal C}^\star(\iota)$ under the identification described by the vertical equalities. Note that by Theorem \ref{TH:Env_C-functor-v-AugSteAlg} the last three horizontal arrows are morphisms of augmented stereotype algebras.

Further we shall need two lemmas.

\blm\label{LM:idealy-v-C(M,B(X))}\footnote{Lemma \ref{LM:idealy-v-C(M,B(X))} was suggested to the author by Robert Israel.} Let $M$ be a paracompact locally compact topological space, and $X$ a finite dimensional vector space over $\C$. Each closed two-sided ideal  $I$ in the algebra ${\mathcal C}\big(M,{\mathcal B}(X)\big)$ has the form
$$
I=\{f\in {\mathcal C}\big(M,{\mathcal B}(X)\big):\ \forall t\in N\ f(t)=0\}
$$
where $N$ is a closed subset in $M$, and the corresponding quotient algebra has the form
$$
{\mathcal C}\big(M,{\mathcal B}(X)\big)/I={\mathcal C}\big(N,{\mathcal B}(X)\big).
$$
\elm
\bpr
For each point $t\in M$ we put $I(t)=\{f(t);\ f\in I\}$ and
$$
N=\{t\in M:\ I(t)=0\},\qquad I(N)=\{f\in {\mathcal C}\big(M,{\mathcal B}(X)\big):\ f\big|_N=0\}.
$$
Obviously, $I(N)$ is a closed two-sided ideal in ${\mathcal C}\big(M,{\mathcal B}(X)\big)$, and $I\subseteq I(N)$. Let us show that the inverse inclusion is also true. Take $f\in I(N)$. For each point $t\notin N$ we have $I(t)\ne 0$, i.e. $I(t)$ is a non-zero two sided ideal in ${\mathcal B}(X)$. Since ${\mathcal B}(X)$ is a simple algebra, $I(t)={\mathcal B}(X)$. Hence
$$
\forall t\notin N\quad \exists g_t\in I(t)={\mathcal B}(X)\quad g_t(t)=f(t).
$$
Now for each compact set $T\subseteq M$ and for any $\e>0$ we can find a partition of unity $\{\eta_t;\ t\in T\}$ on $T$ such that
$$
\sup_{s\in T}\norm{f(s)-\sum_{t\in T}\eta_t(s)\cdot g_t(s)}<\e.
$$
The sums $\sum_{t\in T}\eta_t\cdot g_t$ belong to $I$, and since $I$ is a closed ideal, we have: $f\in I$.
\epr

\blm\label{LM:prod-A/prod-I}
Let $\{A_\sigma;\sigma\in S\}$ be a family of stereotype algebras,
$$
A=\prod_{\sigma\in S}A_\sigma,
$$
its product, $I$ a closed two sided ideal in $A$, and
$$
I_\sigma=\{x_\sigma;\ x\in I\}
$$
its projections at the components $A_\sigma$. Then
\bit{

\item[(i)] each $I_\sigma$ is a closed two sided ideal in $A_\sigma$,

\item[(ii)] $I$ can be recovered from $I_\sigma$ by the formula
\beq\label{I=prod-I_s}
I=\prod_{\sigma\in S}I_\sigma,
\eeq

\item[(iii)] there exists a unique isomorphism
$$
A/I\cong\prod_{\sigma\in S}(A_\sigma/I_\sigma).
$$
}\eit
\elm
\bpr
For each $\sigma\in S$ we define the embedding $\iota_\sigma:A_\sigma\to A$
$$
\iota_\sigma(p)_\tau=\begin{cases}p, & \tau=\sigma\\ 0, & \tau\ne\sigma \end{cases}.
$$
It is mupltiplicative, but not unital, since it turns the unit $1_\sigma\in A_\sigma$ not into the unit $1\in A$, but into the family
$$
\iota_\sigma(1_\sigma)_\tau=\begin{cases}1_\sigma, & \tau=\sigma\\ 0, & \tau\ne\sigma \end{cases}.
$$
Let us note that
\beq\label{I_sigma=iota_sigma^(-1)(I)}
I_\sigma=\iota_\sigma^{-1}(I).
\eeq
Indeed,
\begin{eqnarray*}
p\in I_\sigma\quad\Longleftrightarrow\quad \exists x\in I\quad p=x_\sigma
\quad\Longleftrightarrow\quad \exists x\in I\quad  \iota_\sigma(p)=x\cdot\iota_\sigma(1_\sigma)
\quad\Longleftrightarrow\\ \Longleftrightarrow\quad \iota_\sigma(p)\in I\quad\Longleftrightarrow\quad p\in \iota_\sigma^{-1}(I).
\end{eqnarray*}

1. From \eqref{I_sigma=iota_sigma^(-1)(I)} it follows immediately that $I_\sigma$ is a closed subset in $A_\sigma$. On the other hand, for each $p\in I_\sigma$ and $q\in A_\sigma$ we have
$$
\iota_\sigma(p\cdot q)=\underbrace{\iota_\sigma(p)}_{\scriptsize\begin{matrix}\text{\rotatebox{90}{$\owns$}} \\ I\end{matrix}}\cdot \underbrace{\iota_\sigma(q)}_{\scriptsize\begin{matrix}\text{\rotatebox{90}{$\owns$}} \\ A\end{matrix}}\in I
\quad\overset{\eqref{I_sigma=iota_sigma^(-1)(I)}}{\Longrightarrow}\quad p\cdot q\in I_\sigma.
$$
And similarly, $q\cdot p\in I_\sigma$. This proves (i).

2. Formula \eqref{I=prod-I_s} is obvious.

3. The formula
$$
\varPhi\{a_\sigma+I_\sigma;\ \sigma\in S\}=\{a_\sigma;\ \sigma\in S\}+I
$$
correctly defines a mapping
$$
\varPhi:\prod_{\sigma\in S}(A_\sigma/I_\sigma)\to A/I
$$
and it is easy to see that this is an isomorphism of locally convex spaces.
\epr

\blm\label{LM:stroenie-Coker(xi)}
The cokernel of the mapping $\xi$ in \eqref{DEF:xi:Env-C*(H)->Env-C*(Z-times-K)} in the category $\AugInvSteAlg$ of augmented involutive stereotype algebras has the form
\beq\label{stroenie-Coker(xi)}
\Coker(\xi)\cong
\prod_{\sigma\in\widehat{K}}{\mathcal C}\Big(M_\sigma,{\mathcal B}(X_\sigma)\Big)
\eeq
where $\{M_\sigma,\sigma\in\widehat{K}\}$ is a family of closed subsets in  $\widehat{Z}$.
\elm
\bpr By Theorem \ref{PROP:AugSteAlg-imeet-ker} $\Coker(\xi)$ is a quotient algebra of the algebra $\prod_{\sigma\in\widehat{K}}{\mathcal C}\big(\widehat{Z},{\mathcal B}(X_\sigma)\big)$ by some closed two-sided ideal $I$. By Lemma \ref{LM:prod-A/prod-I} such a quotient algebra is isomorphic to the product of the quotient algebras ${\mathcal C}\big(\widehat{Z},{\mathcal B}(X_\sigma)\big)$ by the ideals $I_\sigma$. And by Lemma  \ref{LM:idealy-v-C(M,B(X))} these quotient algebras are isomorphic to the algebras ${\mathcal C}\big(M_\sigma,{\mathcal B}(X_\sigma)\big)$:
\begin{eqnarray*}
\Coker(\xi)\cong \bigg(\Big(\prod_{\sigma\in\widehat{K}}{\mathcal C}\big(\widehat{Z},{\mathcal B}(X_\sigma)\big)\Big)/I\bigg)^\vartriangle\cong \bigg(\prod_{\sigma\in\widehat{K}}\Big({\mathcal C}\big(\widehat{Z},{\mathcal B}(X_\sigma)\big)/I_\sigma\Big)\bigg)^\vartriangle\cong\\ \cong
\bigg(\prod_{\sigma\in\widehat{K}}{\mathcal C}\big(M_\sigma,{\mathcal B}(X_\sigma)\big)\bigg)^\vartriangle
\cong
\prod_{\sigma\in\widehat{K}}{\mathcal C}\big(M_\sigma,{\mathcal B}(X_\sigma)\big).
\end{eqnarray*}
\epr

Let $Z\cdot K$ be a buildup of an abelian locally compact group $Z$ with the help of a compact group $K$. Consider the chain of homomorphisms
\beq\label{H->Z-times-K->Z-cdot-K}
 \xymatrix  
{
H\ar[r]^{\iota} & Z\times K\ar[r]^{\varkappa} & Z\cdot K
}
\eeq
where $H$ and $\iota$ are defined in Lemma \ref{LM:Z-cdot-K=(Z-times-K)/C}, and
$$
\varkappa=\coker\iota.
$$
(by Lemma \ref{LM:Z-cdot-K=(Z-times-K)/C} $\varkappa$ can be treated as a mapping from $Z\times K$ into $Z\cdot K$). By Theorem \ref{TH:coker-groups=>coker-group-algebras} in the corresponding chain of morphisms of group algebras
\beq\label{C*(H)->C*(Z-times-K)->C*(Z-cdot-K)}
 \xymatrix  @R=2.pc @C=4.pc
{
{\mathcal C}^\star(H)\ar[r]^{\iota^*} & {\mathcal C}^\star(Z\times K)\ar[r]^{\varkappa^*} & {\mathcal C}^\star(Z\cdot K)
}
\eeq
the second morphism is a cokernel of the first one:
$$
\varkappa^*=\coker\iota^*.
$$
In other words,
$$
{\mathcal C}^\star(Z\cdot K)=\Coker\iota^*.
$$

Consider the chain obtained from \eqref{C*(H)->C*(Z-times-K)->C*(Z-cdot-K)} by  applying the functor of continuous envelope:
\beq\label{Env-C*(H)->Env-C*(Z-times-K)->Env-C*(Z-cdot-K)}
 \xymatrix  @R=2.pc @C=4.pc
{
\Env_{\mathcal C} {\mathcal C}^\star(H)\ar[r]^{\Env_{\mathcal C} (\iota^*)} & \Env_{\mathcal C} {\mathcal C}^\star(Z\times K)\ar[r]^{\Env_{\mathcal C} (\varkappa^*)} & \Env_{\mathcal C} {\mathcal C}^\star(Z\cdot K)
}
\eeq
By Theorem \ref{TH:Env_C-functor-v-AugSteAlg} this is a chain of morphisms of augmented involutive stereotype algebras.

\blm\label{LM:Coker(Env_C(iota^*))->Env_C-C^*(Z-cdot-K)}
The mapping
$$
\widehat{\varkappa}:\Coker(\Env_{\mathcal C} (\iota^*))\to \Env_{\mathcal C} {\mathcal C}^\star(Z\cdot K)
$$
is an isomorphism of stereotype algebras.
\elm
\bpr
1. We first prove that the mapping $\widehat{\varkappa}$ is open. From Lemma \ref{LM:stroenie-Coker(xi)} it follows that the topology of the algebra $\Coker(\Env_{\mathcal C} (\iota^*))$ is generated (without pseudosaturation) by the inverse images of the open sets under the morphisms (continuous homomorphisms without taking into account the augmentation) $\ph:\Coker(\Env_{\mathcal C} (\iota^*))\to B$ into different $C^*$-algebras. Hence for proving the openness of the mapping $\Coker(\Env_{\mathcal C} (\iota^*))\to \Env_{\mathcal C} {\mathcal C}^\star(Z\cdot K)$ it is sufficient to show that each morphism $\ph:\Coker(\Env_{\mathcal C} (\iota^*))\to B$ into an arbitrary $C^*$-algebra can be extended to some morphism $\ph':\Env_{\mathcal C} {\mathcal C}^\star(Z\cdot K)\to B$.

Take such a morphism $\ph:\Coker(\Env_{\mathcal C} (\iota^*))\to B$ and consider the diagram
\beq\label{ph:Coker(Env_C(iota^*))->B}
 \xymatrix  @R=2.pc @C=3.pc
{
H\ar[d]_{\delta_H}\ar[r]^{\iota} & Z\times K\ar[d]_{\delta_{Z\times K}}\ar[rrr]^{\varkappa} & & & Z\cdot K\ar[dl]_{\delta_{Z\cdot K}}\ar@{-->}@/^2ex/[ddd]^{\psi} \\
{\mathcal C}^\star(H)\ar[d]_{\env_{\mathcal C}{\mathcal C}^\star(H)}\ar[r]^{\iota^*} & {\mathcal C}^\star(Z\times K)\ar[d]_{\env_{\mathcal C}{\mathcal C}^\star(Z\times K)}\ar[r]^{\coker\iota^*} & \Coker\iota^*\ar@{=}[r]^{\eqref{coker-C^star(lambda)=C^star(varkappa)}} \ar[d]_{\Coker\frac{\Env_{\mathcal C}\iota^*,\env_{\mathcal C}{\mathcal C}^\star(H)}{\env_{\mathcal C}{\mathcal C}^\star(Z\times K),\iota^*}} & {\mathcal C}^\star(Z\cdot K)
\ar@{-->}@/^4ex/[rdd]^{\psi'}\ar[d]_{\env_{\mathcal C} {\mathcal C}^\star(Z\cdot K)} &
\\
\Env_{\mathcal C}{\mathcal C}^\star(H)\ar[r]^{\Env_{\mathcal C}\iota^*} & \Env_{\mathcal C}{\mathcal C}^\star(Z\times K)\ar[r]^{\coker\Env_{\mathcal C}\iota^*} & \Coker\Env_{\mathcal C}\iota^*\ar@/_6ex/[drr]_{\ph}\ar[r]^{\widehat{\varkappa}} &
\Env_{\mathcal C} {\mathcal C}^\star(Z\cdot K)\ar@{-->}@/_2ex/[dr]_{\ph'}\\
&& & & B
}
\eeq
For an element $t\in H$ we have 
$$
\ph\bigg(\Coker\frac{\Env_{\mathcal C}\iota^*,\env_{\mathcal C}{\mathcal C}^\star(H)}{\env_{\mathcal C}{\mathcal C}^\star(Z\times K),\iota^*}\Big(\coker\iota^*(\delta_{Z\times K}^{\iota(t)})\Big)\bigg)=
\ph\bigg(\Coker\frac{\Env_{\mathcal C}\iota^*,\env_{\mathcal C}{\mathcal C}^\star(H)}{\env_{\mathcal C}{\mathcal C}^\star(Z\times K),\iota^*}\Big(\underbrace{\delta_{Z\cdot K}^{\varkappa\big(\iota(t)\big)}}_{\scriptsize\begin{matrix}\|\\ 1 \end{matrix}}\Big)\bigg)=\ph(1)=1,
$$
so we can conclude that the representation of the group $Z\times K$ in the algebra $B$
$$
\ph\circ\Coker\frac{\Env_{\mathcal C}\iota^*,\env_{\mathcal C}{\mathcal C}^\star(H)}{\env_{\mathcal C}{\mathcal C}^\star(Z\times K),\iota^*}\circ\coker\iota^*\circ\delta_{Z\times K}:Z\times K\to B
$$
can be extended to some representation $\psi:Z\cdot K\to B$. By the main property of group algebras \cite[Theorem 10.12]{Akbarov} this represetation generates a morphism of stereotype algebras $\psi':{\mathcal C}^\star(Z\cdot K)\to B$. Since $B$ is a  $C^*$-algebra, the morphism $\psi'$ can be (uniquely) extended to the envelope as a morphism $\ph':\Env_{\mathcal C}{\mathcal C}^\star(Z\cdot K)\to B$. The whole diagram remains commutative since $\coker\iota^*$, $\env_{\mathcal C} {\mathcal C}^\star(Z\cdot K)$, $\widehat{\varkappa}$ are epimorphisms, and the image $\delta_{Z\times K}$ has dense linear span in ${\mathcal C}^\star(Z\times K)$.

2. Next we prove that the mapping $\widehat{\varkappa}$ is injective. Take $x\in \Coker\Env_{\mathcal C}\iota^*$, such that $x\ne 0$. By Lemma  \ref{LM:stroenie-Coker(xi)} we can identify $\Coker\Env_{\mathcal C}\iota^*$ with the algebra of the form  $\prod_{\sigma\in\widehat{K}}{\mathcal C}(M_\sigma,{\mathcal B}(X_\sigma))$. Hence the inequality $x\ne 0$ in $\Coker\Env_{\mathcal C}\iota^*$ implies the existence of a morphism $\ph:\Coker\Env_{\mathcal C}\iota^*\to B$ into a $C^*$-algebra $B={\mathcal B}(X_\sigma)$ such that $\ph(x)\ne 0$.
\beq\label{Coker(Env_C(iota^*))->Env_C-C^*(Z-cdot-K)-1}
 \xymatrix  @R=2.5pc @C=1.5pc
{
&
{\mathcal C}^\star(Z\cdot K)\ar@/_4ex/[ddl]_{\Coker\frac{\Env_{\mathcal C}\iota^*,\env_{\mathcal C}{\mathcal C}^\star(H)}{\env_{\mathcal C}{\mathcal C}^\star(Z\times K),\iota^*}}\ar@/^4ex/[ddr]^{\env_{\mathcal C}{\mathcal C}^\star(Z\cdot K)}\ar@{-->}[d]_{\ph'} & \\
& B & \\
\Coker\Env_{\mathcal C}\iota^*\ar@/_4ex/[rr]_{\widehat{\varkappa}}\ar[ur]_{\ph} & & \Env_{\mathcal C}{\mathcal C}^\star(Z\cdot K)\ar@{-->}[ul]_{\ph''}
}
\eeq
Consider the composition
$$
\ph'=\ph\circ\Coker\frac{\Env_{\mathcal C}\iota^*,\env_{\mathcal C}{\mathcal C}^\star(H)}{\env_{\mathcal C}{\mathcal C}^\star(Z\times K),\iota^*}.
$$
This is a morphism into the $C^*$-algebra $B$, therefore it can be (uniquely) extended to a morphism $\ph''$ at the envelope.

Note further that in Diagram \eqref{ph:Coker(Env_C(iota^*))->B} the morphism $\env_{\mathcal C}{\mathcal C}^\star(Z\times K)$ is dense (being an envelope), and the morphism $\coker\Env_{\mathcal C}\iota^*$ is dense by the construction of the cokernel in Theorem \ref{PROP:AugSteAlg-imeet-ker}. Hence in the equality
$$
\Coker\frac{\Env_{\mathcal C}\iota^*,\env_{\mathcal C}{\mathcal C}^\star(H)}{\env_{\mathcal C}{\mathcal C}^\star(Z\times K),\iota^*}\circ\coker\iota^*=\coker\Env_{\mathcal C}\iota^*\circ\env_{\mathcal C}{\mathcal C}^\star(Z\times K)
$$
the right side consists of dense morphisms. As a corollary, the right side is a dense morphism, and therefore the morphism at the left side, $\Coker\frac{\Env_{\mathcal C}\iota^*,\env_{\mathcal C}{\mathcal C}^\star(H)}{\env_{\mathcal C}{\mathcal C}^\star(Z\times K),\iota^*}$, is dense as well. We can conclude from this that in Diagram \eqref{Coker(Env_C(iota^*))->Env_C-C^*(Z-cdot-K)-1} not only the perimeter and the upper inner triangles are commutative, but the lower inner triangle is commutative as well:
$$
\ph'\circ\widehat{\varkappa}=\ph.
$$
Now we obtain
$$
\ph'(\widehat{\varkappa}(x))=\ph(x)\ne0\quad\Longrightarrow\quad \widehat{\varkappa}(x)\ne 0,
$$
and this is what we need.

3. Now let us prove that the mapping $\widehat{\varkappa}$ is dense. The image of the mapping $\delta:Z\cdot K\to {\mathcal C}^\star(Z\cdot K)$ has a dense linear span in ${\mathcal C}^\star(Z\cdot K)$, hence in the composition with ${\mathcal C}^\star(Z\cdot K)\to\Env_{\mathcal C} {\mathcal C}^\star(Z\cdot K)$ we again have a mapping whose image has a dense linear span in $\Env_{\mathcal C} {\mathcal C}^\star(Z\cdot K)$. As a corollary, its composition with $G\to G/H=Z\cdot K$,
$$
\xymatrix  @R=2.pc @C=4.pc
{
G\ar[r]& G/H=Z\cdot K\ar[r] & \Env_{\mathcal C} {\mathcal C}^\star(Z\cdot K)
}
$$
yields a map whose image has a dense linear span in $\Env_{\mathcal C} {\mathcal C}^\star(Z\cdot K)$. But this map can be represented as a composition
$$
\xymatrix  @R=2.pc @C=4.pc
{
G\ar[r]& \Coker(\Env_{\mathcal C} (\iota^*))\ar[r] & \Env_{\mathcal C} {\mathcal C}^\star(Z\cdot K)
}
$$
and we see that the span of this mapping's image needs to be dense in $\Env_{\mathcal C} {\mathcal C}^\star(Z\cdot K)$. Hence the mapping
$$
\xymatrix  @R=2.pc @C=4.pc
{
\Coker(\Env_{\mathcal C} (\iota^*))\ar[r] & \Env_{\mathcal C} {\mathcal C}^\star(Z\cdot K)
}
$$
possesses the same property.

4. We saw that the mapping of stereotype spaces
$$
\xymatrix  @R=2.pc @C=4.pc
{
\Coker(\Env_{\mathcal C} (\iota^*))\ar[r] & \Env_{\mathcal C} {\mathcal C}^\star(Z\cdot K)
}
$$
is injective, open and dense. On the other hand, Lemma \ref{LM:stroenie-Coker(xi)} implies that the domain of this mapping is a complete locally convex space. Thus, this mapping must be an isomorphism. Therefore the domain of this mapping is isomorphic to its range.
\epr

From Lemmas \ref{LM:Coker(Env_C(iota^*))->Env_C-C^*(Z-cdot-K)} and \ref{LM:stroenie-Coker(xi)} we have

\btm\label{TH:Env_C-C*(Z-cdot-K)}
Let $Z\cdot K$ be a buildup of an abelian locally compact group $Z$ with the help of a compact group $K$. Then the continuous envelope $\Env_{\mathcal C} {\mathcal C}^\star(Z\cdot K)$ of its group algebra has the form
\beq\label{Env_C-C*(Z-cdot-K)}
\Env_{\mathcal C} {\mathcal C}^\star(Z\cdot K)\cong
\prod_{\sigma\in\widehat{K}}{\mathcal C}\Big(M_\sigma,{\mathcal B}(X_\sigma)\Big),
\eeq
where $\{M_\sigma;\ \sigma\in\widehat{K}\}$ is a family of closed subsets in the Pontryagin dual group $\widehat{Z}$ of the group $Z$.
\etm

The following two propositions are special cases of Proposition 5.31 and Lemma 5.53 in \cite{Akbarov-C^infty-2}. The proof of these statements contained errors in  \cite{Akbarov-C^infty-2}\footnote{See \cite[Errata]{Akbarov-C^infty-2}.}, and we repair these errors in the special case which is interesting for us here.

\bcor\label{COR:nepr-Env-SIN-gruppy=LCS-lim} Let $G=Z\cdot K$ be a buildup of an abelian locally compact group $Z$ with the help of a compact group $K$. Then the continuous envelope of the group algebra ${\mathcal C}^\star(G)$ coincides with the locally compact Kuznetsova envelope, i.e. with the projective limit of its $C^*$-quotient algebras in the category of locally convex spaces (and in the category of topological algebras):
\beq\label{E(C*(G))=LCS-leftlim}
\Env_{\mathcal C}{\mathcal C}^\star(G)={\tt LCS}\text{-}\kern-3pt\projlim_{p\in{\mathcal P}({\mathcal C}^\star(G))}{\mathcal C}^\star(G)/p
\eeq
\ecor
\bpr
This is seen from \eqref{Env_C-C*(Z-cdot-K)}.
\epr

\bcor\label{COR:o-diagonali-beta}
Let $G=Z\cdot K$ be a buildup of an abelian locally compact group $Z$ with the help of a compact group $K$. Then the diagonal $\beta$ of the diagram
$$
\xymatrix @R=6.pc @C=8.pc 
{
{\mathcal C}^\star(G)\circledast {\mathcal C}^\star(G)\ar@{-->}[dr]^{\beta}
\ar[d]_{@}\ar[r]^(.4){\env_{\mathcal C} {{\mathcal C}^\star(G)}\circledast \env_{\mathcal C} {{\mathcal C}^\star(G)}} & \Env_{\mathcal C} {\mathcal C}^\star(G)\circledast \Env_{\mathcal C} {\mathcal C}^\star(G)\ar[d]_{@} \\
{\mathcal C}^\star(G)\odot {\mathcal C}^\star(G)\ar[r]^(.4){\env_{\mathcal C} {{\mathcal C}^\star(G)}\odot \env_{\mathcal C} {{\mathcal C}^\star(G)}} & \Env_{\mathcal C} {\mathcal C}^\star(G)\odot \Env_{\mathcal C} {\mathcal C}^\star(G)
}
$$
is a dense epimorphism.
\ecor
\bpr
From \eqref{Env_C-C*(Z-cdot-K)} we have that since the image of ${\mathcal C}^\star(Z\cdot K)$ under the mapping $\env_{\mathcal C}{\mathcal C}^\star(Z\cdot K)$ is dense in the space $\prod_{\sigma\in\widehat{K}}{\mathcal C}\big(M_\sigma,{\mathcal B}(X_\sigma)\big)$, under the projection on each finite product
$\prod_{\sigma\in S}{\mathcal C}\big(M_\sigma,{\mathcal B}(X_\sigma)\big)$ (where  $S\subseteq\widehat{K}$ is an arbitrary finite set), the mapping
$$
{\mathcal C}^\star(Z\cdot K)\to \prod_{\sigma\in S}{\mathcal C}\big(M_\sigma,{\mathcal B}(X_\sigma)\big)
$$
is dense as well. This implies that for each finite sets $S,T\subseteq\widehat{K}$ the arising mapping of the tensor products
\begin{eqnarray*}
{\mathcal C}^\star(Z\cdot K)\circledast{\mathcal C}^\star(Z\cdot K)\to \prod_{\sigma\in S}{\mathcal C}\big(M_\sigma,{\mathcal B}(X_\sigma)\big)\circledast
\prod_{\tau\in S}{\mathcal C}\big(M_\tau,{\mathcal B}(X_\tau)\big)=\\=\prod_{\sigma\in S,\tau\in T}{\mathcal C}\big(M_\sigma,{\mathcal B}(X_\sigma)\big)\circledast {\mathcal C}\big(M_\tau,{\mathcal B}(X_\tau)\big)
\end{eqnarray*}
is again dense.

Note further that ${\mathcal C}\big(M_\sigma,{\mathcal B}(X_\sigma)\big)$ and ${\mathcal C}\big(M_\tau,{\mathcal B}(X_\tau)\big)$ are Fr\'echet spaces with the classical approximation property. Hence by \cite[Theorem 7.21]{Akbarov} their stereotype injective tensor product coincides with their classical injective tensor product:
$$
{\mathcal C}\big(M_\sigma,{\mathcal B}(X_\sigma)\big)\odot {\mathcal C}\big(M_\tau,{\mathcal B}(X_\tau)\big)=
{\mathcal C}\big(M_\sigma,{\mathcal B}(X_\sigma)\big)\check{\otimes} {\mathcal C}\big(M_\tau,{\mathcal B}(X_\tau)\big).
$$
At the same time the algebraic tensor product ${\mathcal C}\big(M_\sigma,{\mathcal B}(X_\sigma)\big)\otimes{\mathcal C}\big(M_\tau,{\mathcal B}(X_\tau)\big)$ is dense in \break ${\mathcal C}\big(M_\sigma,{\mathcal B}(X_\sigma)\big)\check{\otimes} {\mathcal C}\big(M_\tau,{\mathcal B}(X_\tau)\big)$. This implies that the projective stereotype tensor product
${\mathcal C}\big(M_\sigma,{\mathcal B}(X_\sigma)\big)\circledast {\mathcal C}\big(M_\tau,{\mathcal B}(X_\tau)\big)$ is also dense in ${\mathcal C}\big(M_\sigma,{\mathcal B}(X_\sigma)\big)\check{\otimes} {\mathcal C}\big(M_\tau,{\mathcal B}(X_\tau)\big)$. This means in its turn that the composition of the mappings
$$
{\mathcal C}^\star(Z\cdot K)\circledast{\mathcal C}^\star(Z\cdot K)\to \prod_{\sigma\in S,\tau\in T}{\mathcal C}\big(M_\sigma,{\mathcal B}(X_\sigma)\big)\circledast {\mathcal C}\big(M_\tau,{\mathcal B}(X_\tau)\big)\to \prod_{\sigma\in S,\tau\in T}{\mathcal C}\big(M_\sigma,{\mathcal B}(X_\sigma)\big)\odot {\mathcal C}\big(M_\tau,{\mathcal B}(X_\tau)\big)
$$
is again a dense mapping.

And this is true for each finite sets $S,T\subseteq\widehat{K}$. We can conclude that the mapping into the infinite product
\begin{eqnarray*}
{\mathcal C}^\star(Z\cdot K)\circledast{\mathcal C}^\star(Z\cdot K)\to \prod_{\sigma,\tau\in \widehat{K}}{\mathcal C}\big(M_\sigma,{\mathcal B}(X_\sigma)\big)\odot {\mathcal C}\big(M_\tau,{\mathcal B}(X_\tau)\big)
\stackrel{\text{\cite[(7.21)]{Akbarov}}}{=}
\\=
\prod_{\sigma\in \widehat{K}}{\mathcal C}\big(M_\sigma,{\mathcal B}(X_\sigma)\big)\odot \prod_{\tau\in \widehat{K}}{\mathcal C}\big(M_\tau,{\mathcal B}(X_\tau)\big)=
\Env_{\mathcal C} {\mathcal C}^\star(Z\cdot K)\odot\Env_{\mathcal C} {\mathcal C}^\star(Z\cdot K)
\end{eqnarray*}
is again dense.
\epr

\btm\label{TH:Env_C-C^*(Z-cdot-K)-odot-Hopf} Let $G=Z\cdot K$ be a buildup of an abelian locally compact group $Z$ with the help of a compact group $K$. Then the continuous envelope $\Env_{\mathcal C} {\mathcal C}^\star(G)$ of its group algebra ${\mathcal C}^\star(G)$ is an involutive Hopf algebra in the category of stereotype spaces $(\tt{Ste},\odot)$.
\etm
\bpr
The proof of this statement repeats almost without changes the proof of Theorem 5.51 in \cite{Akbarov-C^infty-2}. The only difference is that Corollaries  \ref{COR:nepr-Env-SIN-gruppy=LCS-lim} and \ref{COR:o-diagonali-beta} replace in this reasoning Proposition 5.31 (with Formula (5.65)) and Lemma 5.53 from \cite{Akbarov-C^infty-2}.
\epr

\section{Acknowledgements.} The author thanks Robert Israel for a useful consultation.

\end{document}